\documentclass[12pt,reqno]{amsart}
\usepackage{amsmath,amsthm,amssymb,amsfonts,amscd}
\usepackage{mathrsfs}
\usepackage{bbm}
\usepackage{booktabs}
\usepackage{float}
\usepackage[toc,page]{appendix}
\usepackage{makecell}
\usepackage{bbding}
\usepackage[backref=page]{hyperref}
\usepackage{geometry}
\usepackage{color}
\usepackage{xcolor}
\usepackage{picture,epic}
\usepackage{tikz}

\usepackage{etoolbox}
\usepackage{orcidlink}

\numberwithin{equation}{section}

\theoremstyle{plain}
\newtheorem{theorem}{Theorem}[section]
\newtheorem{lemma}[theorem]{Lemma}

\newtheorem{proposition}[theorem]{Proposition}

\theoremstyle{definition}

\theoremstyle{remark}
\newtheorem{remark}[theorem]{Remark}

\renewcommand{\Re}{\operatorname{Re}}

\newcommand{\rank}{\operatorname{rank}}

\newcommand{\dd}{\mathrm{d}}

   \DeclareFontFamily{U}{wncy}{}
    \DeclareFontShape{U}{wncy}{m}{n}{<->wncyr10}{}
    \DeclareSymbolFont{mcy}{U}{wncy}{m}{n}
    \DeclareMathSymbol{\Sh}{\mathord}{mcy}{"58}

\makeatletter
\def\@tocline#1#2#3#4#5#6#7{\relax
  \ifnum #1>\c@tocdepth 
  \else
    \par \addpenalty\@secpenalty\addvspace{#2}%
    \begingroup \hyphenpenalty\@M
    \@ifempty{#4}{%
      \@tempdima\csname r@tocindent\number#1\endcsname\relax
    }{
      \@tempdima#4\relax
    }
    \parindent\z@ \leftskip#3\relax \advance\leftskip\@tempdima\relax
    \rightskip\@pnumwidth plus4em \parfillskip-\@pnumwidth
    #5\leavevmode\hskip-\@tempdima
      \ifcase #1
       \or\or \hskip 1em \or \hskip 2em \else \hskip 3em \fi%
      #6\nobreak\relax
    \hfill\hbox to\@pnumwidth{\@tocpagenum{#7}}\par
    \nobreak
    \endgroup
  \fi}
\makeatother

\begin{document}

\title[$L$-derivatives]{$L$-derivatives of the fixed elliptic curve over rank-one imaginary quadratic fields}

\author{Shenghao Hua~\orcidlink{0000-0002-7210-2650}}
\email{huashenghao@vip.qq.com}


\begin{abstract}
There is a one-sided central limit theorem for the logarithms of $L$-derivatives of a fixed rational non-CM elliptic curve $E$ over imaginary quadratic fields of rank one, analogous to a result of Radziwi\l\l\ and Soundararajan.
There are also many $L$-derivatives that are not small.
\end{abstract}

\keywords{elliptic curve, imaginary quadratic field, $L$-function}

\subjclass[2020]{11G05,11G50,14H15,11R42}

\maketitle

\section{Introduction} \label{sec:Intr}
Let \( E \) be an elliptic curve defined over the rational numbers \(\mathbb{Q}\), and let \( K \) be a number field containing \(\mathbb{Q}\). The base change of \( E \) to \( K \) allows us to consider the group of \( K \)-rational points \( E(K) \).
When \( K=\mathbb{Q}(\sqrt{d}) \) is an imaginary quadratic field with square-free integer $d<0$, the base change of $E$ from \(\mathbb{Q}\) to \( K \) lifts its classical modular form to a Bianchi modular form over \( \mathbb{Q}(\sqrt{d}) \) via Langlands base change theory. These base change forms have a clear structural role within Bianchi modular forms, helping to better understand modularity and Galois representations of elliptic curves over \( \mathbb{Q}(\sqrt{d}) \), and linking classical modular forms to the Langlands program.
On the other hand, extending the base field to \( K \) typically increases the set of \( K \)-rational points on $E$, and arithmetic invariants such as the fundamental units of the Mordell--Weil group may change, thereby enriching the curve's arithmetic properties.

Here, we take the central point to be $s=1$.
For $\Re s>\frac{3}{2}$, the normalized Hasse--Weil $L$-function of $E$ over $\mathbb{Q}$ is
\[
L(s,E) =
L(s,E/\mathbb{Q})
:=\prod_{p \nmid N} \left(1 - a_p p^{-s} + p^{1 - 2s} \right)^{-1} \prod_{p \mid N} L_p(p^{-s})^{-1}
=\sum_{n=1}^\infty \frac{a_n}{n^s},
\]
where \( a_p = p + 1 - \#E(\mathbb{F}_p) \) for primes \( p \) of good reduction, and the local factors \( L_p(p^{-s}) \) at bad primes are determined by the type of reduction of \( E \) at \( p \).
From Hasse's bound, we know that
\[
|a_p| \leq 2 \sqrt{p}.
\]
The completed \(L\)-function of \( E/\mathbb{Q} \) is defined by
\[
\Lambda(s,E)
=
N^{s/2} (2\pi)^{-s} \Gamma(s) L(s,E),
\]
where \( N \) is the conductor of \( E \). Thanks to the modularity theorem proved by Wiles and Taylor, with crucial contributions from Breuil, Conrad, and Diamond, \( \Lambda(E, s) \) admits analytic continuation to the entire complex plane and satisfies the functional equation
\[
\Lambda(s,E) = w \, \Lambda(2-s,E),
\]
with \( w = \pm 1 \) the root number.

Any real primitive character to the modulus $n$ must be of the form $\chi_n(\nu)=(\frac{n}{\nu})$ where $n$ is a fundamental discriminant \cite[Theorem 9.13]{MontgomeryVaughan2007}, i.e., a product of pairwise coprime integers of the form $-4$, $\pm 8$, $(-1)^{\frac{p-1}{2}}p$ where $p$ is an odd prime.
Given a fundamental discriminant \( d \), the quadratic twist \( E^{(d)} \) over $\mathbb{Q}$ has Hasse--Weil \( L \)-function
\[
L(s,E^{(d)})
=
\sum_{n=1}^\infty \frac{a_n \left(\frac{d}{n}\right)}{n^s}
=\sum_{n=1}^\infty \frac{a_n \chi_d(n)}{n^s}.
\]
The completed \(L\)-function of \(E^{(d)}\) is defined by
\[
\Lambda(s,E^{(d)}) = N_d^{s/2} (2\pi)^{-s} \Gamma(s) L(s,E^{(d)}),
\]
where \(N_d\) is the conductor of the twist \(E^{(d)}\).
When $(d,2N)=1$, \(N_d=d^2N\).
We also have
\[
\Lambda(s,E^{(d)}) = w_d \, \Lambda(2-s,E^{(d)}),
\]
with root number
\begin{equation}\label{epsilon}
  w_d=w\chi_d(-N).
\end{equation}
For \( d < 0 \) a fundamental discriminant, the base change \(L\)-function satisfies
\[
L(s, E/\mathbb{Q}(\sqrt{d})) = L(s, E) \cdot L(s, E^{(d)}).
\]

For $\epsilon=\pm 1$, let
\begin{equation*}
  \Omega_{\epsilon} := \{d:\textrm{fundamental discriminant }d<0,~ (d,2N)=1,~ w\chi_d(-N)=\epsilon\}.
\end{equation*}
Let $N_0=[8,N]$.
Let \( a \bmod N_0 \) be a residue class, define
\begin{equation*}
\Omega_{\epsilon}(a) := \{d\in \Omega_{\epsilon} : \  d\equiv a\pmod {N_0}, \ d<0\} .
\end{equation*}

To impose a constraint on $\rank E(\mathbb{Q}(\sqrt{d}))$, we consider the quadratic twist $E^{(d)}$, which is another elliptic curve over $\mathbb{Q}$, defined (up to isomorphism) by:
\[
E^{(d)}: y^2 = x^3 + a d^2 x + b d^3,
\]
if $E$ is given by $y^2 = x^3 + ax + b$.

Let Galois group \( G = \mathrm{Gal}(\mathbb{Q}(\sqrt{d})/\mathbb{Q}) = \{1, \sigma\} \).
Then $E(\mathbb{Q}(\sqrt{d}))$ decomposes as a direct sum of two subgroups:
\[
E(\mathbb{Q}(\sqrt{d})) = E^{+}(\mathbb{Q}(\sqrt{d})) \oplus E^{-}(\mathbb{Q}(\sqrt{d})),
\]
where
\[
E^{+}(\mathbb{Q}(\sqrt{d}))
= \{ P \in E(\mathbb{Q}(\sqrt{d})) : \sigma(P) = P \}
=E(\mathbb{Q})
\]
is the subgroup fixed by \( G \), and
\[
E^{+}(\mathbb{Q}(\sqrt{d}))
=\{ P \in E(\mathbb{Q}(\sqrt{d})) : \sigma(P) = -P \}
\cong E^{(d)}(\mathbb{Q}).
\]
Thus we have
\begin{equation*}
  \rank E(\mathbb{Q}(\sqrt{d}))=\rank E(\mathbb{Q})+\rank E^{(d)}(\mathbb{Q}).
\end{equation*}

For an elliptic curve over $\mathbb{Q}$, when the vanishing order of its $L$-function at the central point, i.e., the analytic rank, is at most 1, then the results of Gross--Zagier~\cite{GrossZagier1986} and Kolyvagin~\cite{Kolyvagin1990} imply that the analytic rank coincides with the algebraic rank.
Therefore, for \( i = 0, 1 \), if the analytic rank of \( E \) over \( \mathbb{Q} \) is \( i \), we require that the analytic rank of its quadratic twist \( E^{(d)} \) over \( \mathbb{Q} \) be \( 1 - i \).
Then from
\begin{equation}\label{eqn:L'EK}
  L'(1, E/\mathbb{Q}(\sqrt{d})) = L(1, E) L'(1, E^{(d)}) + L'(1, E) L(1, E^{(d)}),
\end{equation}
we have the following results on $L$-functions.

\begin{theorem}[One-sided central limit theorem]\label{thm:CLT-L}
Let \( E \) be a non-CM elliptic curve over \( \mathbb{Q} \) with conductor \( N \) and root number \( w \), and the vanishing order of the \(L(s,E)\) at the central point is at most 1.
Fix a residue class \( a \bmod N_0 \) be a residue class such that \( a \) is a quadratic residue modulo \( N \), \( a \equiv 1 \) or \( 5 \pmod{8} \),
and for any fundamental discriminant $d<0$ with $d\equiv a \pmod{N_0}$, $ w\chi_d(-N)=\epsilon$.
Let \( V \in \mathbb{R} \).
For large $X$ we have
\[
\left |\left\{
 d \in \Omega_{-w}(a)\mid 20 < -d \leq X,~
\frac{ \log L'(1, E/\mathbb{Q}(\sqrt{d})) - \tfrac{w}{2} \log \log |d| }{ \sqrt{ \log \log |d| } } \geq V
\right\}
\right |
\]
is at most
\[
\left| \{ d \in \Omega_{-w}(a) : -d \leq X \} \right| \cdot \left( \frac{1}{\sqrt{2\pi}} \int_V^{\infty} e^{-x^2/2} \, dx + o(1) \right).
\]
Moreover, when \(L(1,E)\neq 0\), this result also applies to \(L'(1, E^{(d)})\).
\end{theorem}

\begin{remark}
  Here, we mainly employ the method of Radziwi\l\l\ and Soundararajan to establish a one-sided central limit theorem.
  For the large deviations problem, it may be possible to adapt techniques similar to those of Arguin and Bailey~\cite{ArguinBailey2022} and Creighton~\cite{Creighton2025} to control the number of twists with large values of \(\log L'(1, E/\mathbb{Q}(\sqrt{d}))\), say of size comparable to \(\alpha \log \log |d|\), for \(0 <\alpha <2\).
Assuming the Generalized Riemann Hypothesis, a full version of the central limit theorem can be established, following the approach presented in Radziwi\l\l\--Soundararajan~\cite{RadziwillSoundararajan2024}.
\end{remark}

\begin{theorem}\label{thm:lowerbound-L}
Let \( E \) be a non-CM elliptic curve over \( \mathbb{Q} \) with conductor \( N \) and root number \( w \), and the vanishing order of the \(L(s,E)\) at the central point is at most 1.
Fix a residue class \( a \bmod N_0 \) be a residue class such that \( a \) is a quadratic residue modulo \( N \), \( a \equiv 1 \) or \( 5 \pmod{8} \),
and for any fundamental discriminant $d<0$ with $d\equiv a \pmod{N_0}$, $ w\chi_d(-N)=\epsilon$.
For \( -3>d\in  \Omega_{-w}(a)\), there exist positive constants \( c_1', c_2', c_3' \) such that for any small \( \varepsilon > 0 \) and sufficiently large \( X \), we have
\[
\left|\left\{
 d \in \Omega_{1}(a) \mid 20 < -d \leq X,\
L'(1, E/\mathbb{Q}(\sqrt{d})) \geq c_1'
\right\}\right| \gg X^{1 - \varepsilon},
\]
\[
\left|\left\{
 d \in \Omega_{-1}(a) \mid 20 < -d \leq X,\
L'(1, E/\mathbb{Q}(\sqrt{d}))\geq c_2' \log|d| + c_3'
\right\}\right| \gg X^{1 - \varepsilon}.
\]
Moreover, when \(L(1,E)\neq 0\), this result also applies to \(L'(1, E^{(d)})\).
\end{theorem}

\section{Preparation}\label{sec:lemmas}
Let $\Phi(\cdot)\leq 1$ be a smooth non-negative Schwartz class function supported on $[\frac{1}{2},\frac{5}{2}]$ with $\Phi(x)=1$ for $x\in [1,2]$, and for any complex number $s$ let
\begin{equation*}
  \check\Phi(s)=\int_{0}^{\infty}\Phi(x)x^s\dd x.
\end{equation*}
Similarly to Radziwi\l\l--Soundararajan~\cite{RadziwillSoundararajan2015} and Gao--Zhao~\cite{GaoZhao2024}, we have a twisted first moment of central $L$-derivatives for an elliptic curve $E$.

\begin{lemma}\label{lemma:twist}
Let \( u \) be a positive integer such that \( (u\ell, N_0) = 1 \).
For \( w = 1 \), define
\begin{equation*}
  S'(X; u) = \sum_{\substack{d \in \Omega_{-1}(a)}}
  L'\left(1, E^{(d)}\right) \chi_d(u)
  \Phi\left(\frac{-d}{X}\right).
\end{equation*}

Write \( u = u_1 u_2^2 \leq p \) with \( u_1 \) square-free.
Then there exist absolute constants \( C = C(E) \), \( C_1 = C_1(E, \Phi) \), and constants \( C_2(p) \ll 1 \) for all \( p \), such that for any \( \varepsilon > 0 \), we have
\begin{equation*}
  S'(X; u) = \frac{a(u_1)}{4\pi^2 u_1 N_0} \,
  \check{\Phi}(0) X \left( \log \frac{X}{u_1} + C_1 + \sum_{p \mid u} \frac{C_2(p)}{p} \log p \right)
  + O\left(X^{\frac{3}{4} + \varepsilon} u^{1/2}\right).
\end{equation*}

Here, we write \( a(p) = \alpha_p + \beta_p \), where \( \alpha_p \beta_p = p \), and \( |\alpha_p| = |\beta_p| = \sqrt{p} \) for all primes \( p \nmid N \).
\end{lemma}

\begin{proof}
We omit the detailed computation here. Our argument follows the same method as in the proof of Gao--Zhao~\cite[Lemma 3.3]{GaoZhao2024}, which uses normalized notation (with central point \( s = \frac{1}{2} \)) and bounds the Fourier coefficients by the divisor function. Their result is itself a twisted version of Shen~\cite[Theorem 1.2]{Shen2022}. We establish our result by adapting this approach with two modifications.

First, while \cite{GaoZhao2024,Shen2022} consider the full modular group case, the elliptic curve $E$ corresponds to a modular form of level strictly greater than one, but as demonstrated in the proof of Radziwi\l\l--Soundararajan~\cite[Proposition 1]{RadziwillSoundararajan2015}, the strategy of Soundararajan~\cite{Soundararajan2000} also works in this setting.

Second, unlike \cite{GaoZhao2024,Shen2022} which work with fundamental discriminants of the form $8d$ for positive odd square-free $d$, we restrict fundamental discriminants to the form $d \in \Omega_{-1}(a)$. Heath-Brown's large sieve \cite{HeathBrown1995} remains applicable in this setting, but the main term must be adapted to accommodate the new Euler product structure, similar to the adjustment in \cite{RadziwillSoundararajan2015}.

The error term in our case is also $O(X^{3/4+\varepsilon}u^{1/2})$, respect to \cite{GaoZhao2024}.
\end{proof}

\begin{remark}
  Ricotta and Templier previously derived an asymptotic formula with an error term of $O(X^{19/20+\varepsilon})$~\cite{RicottaTemplier2009}, which they used to obtain the average height of Heegner points.
\end{remark}

\begin{lemma}[\cite{RadziwillSoundararajan2015},~Lemma 3]\label{lemma:ortho}
There exists a positive constant \( c' \) such that
\begin{equation*}
  \sum_{p \leq x} \frac{a(p)^2}{p} \log p = x + O\bigl(x \exp(-c' \sqrt{x})\bigr).
\end{equation*}
Furthermore, there exists a constant \( B \) such that
\begin{equation*}
  \sum_{p \leq x} \frac{a(p)^2}{p^2} = \log \log x + B + O\left(\frac{1}{\log x}\right).
\end{equation*}
\end{lemma}

When $w=1$, let \( \mathcal{P} \) denote the set of primes \( p \leq X^{\frac{1}{(\log\log X)^2}} \) such that \( p \nmid N_0 \).
Let \( d \in \Omega_{-1}(a) \) with \( X \leq -d \leq 2X \), and define
\[
\mathcal{P}(d) = \sum_{p \in \mathcal{P}} \frac{a(p)}{p} \chi_d(p).
\]

\begin{lemma}\label{lemma:Pbound}
For any fixed \( V \in \mathbb{R} \) and \( \varepsilon > 0 \), there are \( o(X) \) values of \( d \in \Omega_{-1}(a) \) with \( \frac{X}{\log X} < -d \leq X \) such that
\begin{equation*}
  \mathcal{P}(d) \geq (V - \varepsilon) \sqrt{\log\log X},
\end{equation*}
or
\begin{equation*}
  |\mathcal{P}(d)| \geq \log\log X.
\end{equation*}
\end{lemma}

\begin{proof}
We follow essentially the same proof as in \cite[\S 4]{RadziwillSoundararajan2015}, where the authors work with normalized Fourier coefficients. Their argument based on the Poisson summation formula and Selberg orthogonality, does not rely on the root number of \( E^{(d)} \).
\end{proof}

Let \(\ell\) be a non-negative integer and \(x\) a real number. Define
\[
E_{\ell}(x) = \sum_{j=0}^{\ell} \frac{x^{j}}{j!}.
\]

\begin{lemma}[\cite{RadziwillSoundararajan2015},~Lemma 1]\label{lemma:el}
Let \( \ell \) be a non-negative even integer. Then the function \( E_{\ell}(x) \) is positive and convex on \( \mathbb{R} \). Moreover, for any \( x \leq 0 \), we have
\[
E_{\ell}(x) \geq e^{x}.
\]
Furthermore, if \( \ell \) is a positive even integer and \( x \leq \ell / e^2 \), then
\[
e^{x} \leq \left( 1 + \frac{e^{-\ell}}{16} \right) E_{\ell}(x).
\]
\end{lemma}

\begin{proposition}\label{prop:bound}
With the notation as above, let \( \ell = 20 \lfloor \log \log X \rfloor \), then we have
\begin{equation*}
\sum_{\substack{d \in \Omega_{-1}(a)}}
L'(1, E^{(d)})
E_{\ell}(-\mathcal{P}(d))
\Phi\left(\frac{-d}{X}\right)
\ll X (\log X)^{\frac{1}{2}}\log \log X.
\end{equation*}
\end{proposition}

\section{Proof of the Theorems for \texorpdfstring{$L$}{L}-functions}\label{sec:proof}

\begin{proof}[Proof of Theorem~\ref{thm:CLT-L}]
From \eqref{eqn:L'EK}, when \( w = -1 \),
we apply \cite[Theorem~2]{RadziwillSoundararajan2015}, whose proof is carried out in the setting of residue classes modulo $N_0$.

When $w=1$, we need to establish a similar argument to that of \cite[Theorem~2]{RadziwillSoundararajan2015}.
If \( d \in \Omega_{-1}(a) \) with \( \frac{X}{\log X} < -d \leq X \) satisfies
\[
\log L'(1, E^{(d)}) - \tfrac{1}{2} \log \log X \geq V \sqrt{\log \log X},
\]
then either \( \mathcal{P}(d) \) satisfies the bound in Lemma~\ref{lemma:Pbound}, or
\[
- \log \log X \leq \mathcal{P}(d) \leq (V - \varepsilon ) \sqrt{\log \log X}
\]
but
\[
L'(1, E^{(d)}) (\log X)^{-1/2} \exp(-\mathcal{P}(d)) \geq \exp(\varepsilon \sqrt{\log \log X}).
\]

Recall that \( \ell \) satisfies the condition \( \ell \geq e^2 \left| \mathcal{P}(d) \right| \).
By Lemma~\ref{lemma:el}, we have
\[
L'(1, E^{(d)}) (\log X)^{-1/2} E_\ell (-\mathcal{P}(d)) \gg \exp(\varepsilon \sqrt{\log \log X}).
\]
In light of Proposition~\ref{prop:bound}, this final case occurs with frequency \( o(X) \).
\end{proof}

\begin{proof}[Proof of Theorem~\ref{thm:lowerbound-L}]
From \eqref{eqn:L'EK}, when \( w = -1 \),
applying \cite[Proposition~2]{RadziwillSoundararajan2015}, there exists a constant \( D_1 > 0 \), depending on \( E \), such that
\begin{multline*}
\sum_{d \in \Omega_{1}(a)}
L'(1, E/\mathbb{Q}(\sqrt{d}))
\Phi\left(\frac{-d}{X}\right)
=
L'(1, E)
\sum_{d \in \Omega_{1}(a)} L(1, E^{(d)})
\Phi\left(\frac{-d}{X}\right)
\\= D_1 X + O\left(X^{\frac{7}{8}+\varepsilon}\right).
\end{multline*}

When \( w = 1 \), using Lemma~\ref{lemma:twist}, there exist constants \( D_2, D_3 > 0 \), depending on \( E \), such that
\begin{multline*}
\sum_{d \in \Omega_{-1}(a)}
L'(1, E/\mathbb{Q}(\sqrt{d}))
\Phi\left(\frac{-d}{X}\right)
=
L(1, E)
\sum_{d \in \Omega_{1}(a)}
L'(1, E^{(d)})
\Phi\left(\frac{-d}{X}\right)
\\= D_2 X \log X + D_3 X + O\left(X^{\frac{3}{4}+\varepsilon}\right).
\end{multline*}

Moreover, by Heath-Brown's quadratic large sieve~\cite{HeathBrown1995}, in both cases we have
\[
\sum_{d \in \Omega_{-w}(a)}
L'(1, E/\mathbb{Q}(\sqrt{d}))^2
\Phi\left(\frac{-d}{X}\right)
\ll X^{1 + \frac{\varepsilon}{2}}.
\]

Let \( c_i'< \frac{D_i}{2} \) for \( i = 1, 2, 3 \). Then for \( w = -1 \), we have
\[
\sum_{\substack{d \in \Omega_1(a) \\ L'(1, E/\mathbb{Q}(\sqrt{d})) < c_1'}}
L'(1, E/\mathbb{Q}(\sqrt{d}))
\Phi\left(\frac{-d}{X}\right)
< 2c_1' X,
\]
and thus,
\[
\sum_{\substack{d \in \Omega_1(a) \\ L'(1, E/\mathbb{Q}(\sqrt{d})) \geq c_1'}}
L'(1, E/\mathbb{Q}(\sqrt{d}))
\Phi\left(\frac{-d}{X}\right)
> (D_1 - 2c_1') X.
\]

A similar argument applies for \( w = 1 \).

Now let \( k = \frac{\varepsilon}{3} \). By H\"older's inequality, we obtain
\begin{multline*}
\sum_{\substack{d \in \Omega_1(a) \\ L'(1, E/\mathbb{Q}(\sqrt{d})) \geq c_1'}}
\left|L'(1, E/\mathbb{Q}(\sqrt{d}))
\Phi\left(\frac{-d}{X}\right)\right|^k
\\ \gg
\left(
\sum_{\substack{d \in \Omega_1(a) \\ L'(1, E/\mathbb{Q}(\sqrt{d})) \geq c_1'}}
\left|L'(1, E/\mathbb{Q}(\sqrt{d}))
\Phi\left(\frac{-d}{X}\right)\right|
\right)^{2 - k}
\\ \times \left(
\sum_{\substack{d \in \Omega_1(a) \\ L'(1, E/\mathbb{Q}(\sqrt{d})) \geq c_1'}}
\left|L'(1, E/\mathbb{Q}(\sqrt{d}))
\Phi\left(\frac{-d}{X}\right)\right|^2
\right)^{k - 1}
\gg X^{1 - \frac{\varepsilon}{2}}.
\end{multline*}

Using the convexity bound
\[
L(1, E^{(d)}),~ L'(1, E^{(d)})
\ll d^{\frac{1}{2} + \frac{\varepsilon}{2}},
\]
we deduce
\[
\sum_{\substack{d \in \Omega_1(a) \\ L'(1, E/\mathbb{Q}(\sqrt{d})) \geq c_1'}} 1
\gg X^{1 - \varepsilon}.
\]

A similar conclusion holds for \( w = 1 \). This completes the proof.
\end{proof}

\section{Proof of Proposition~\ref{prop:bound}}
\begin{proof}
Our proof follows a similar strategy to that of \cite[Proposition~5]{RadziwillSoundararajan2015},
but with an additional logarithmic factor, as in \cite{GaoZhao2024}.

For any integer $n$, we write $n = (n)_1 (n)_2^2$,
where $(n)_1$ is square-free.
Let $\tilde{a}(\cdot)$ be a completely multiplicative function defined by $
\tilde{a}(p^\alpha) = a(p)^{\alpha}$ for all $\alpha \geq 0$.
Let $\Omega(n)$ denote the total number of prime factors of $n$, counted with multiplicity.
Let $\omega(\cdot)$ be a multiplicative function defined by $\omega(p^\alpha) = \alpha !$.
Define the function $b(n)$ to be 1 if $n$ has at most $\ell$ prime factors, all of which lie in $\mathcal{P}$, and 0 otherwise.
Then we have
\begin{equation*}
  E_{\ell}(-\mathcal{P}(d))
  =
  \sum_{n}
  \frac{\tilde{a}(n)}{n}
 \frac{(-1)^{\Omega(n)}}{\omega(n)}
 b(n) \chi_d(n).
\end{equation*}

Using Lemma~\ref{lemma:twist}, we obtain
\begin{multline*}
\sum_{\substack{d \in \Omega_{-1}(a)}}
L'(1, E^{(d)})
E_{\ell}(-\mathcal{P}(d))
\Phi\left(\frac{-d}{X}\right)
\\=
\frac{ \check{\Phi}(0) X }{4\pi^2 N_0}
\sum_{n}
\frac{\tilde{a}(n) a((n)_1)}{n (n)_1}
 \frac{(-1)^{\Omega(n)}}{\omega(n)}
 b(n)
\\ \times \left( \log \frac{X}{(n)_1} + C_1 +
\sum_{p \mid n} \frac{C_2(p)}{p} \log p \right)
+ O\left(X^{\frac{3}{4} + 2\varepsilon}\right).
\end{multline*}

We shall concentrate on the terms involving
\(\log\left(\frac{X}{(n)_1}\right)\). The remaining terms, involving
\[
C_1 + \sum_{p \mid n} \frac{C_2(p)}{p} \log p,
\]
can be shown, using the same argument, to contribute
at most \( O\left(X (\log X)^{-\frac{1}{2}}\log\log X\right)\).

Define
\begin{equation*}
  S_1 =
  \sum_{n}
  \frac{\tilde{a}(n) a((n)_1)}{n (n)_1}
\frac{(-1)^{\Omega(n)}}{\omega(n)}
b(n),
\end{equation*}
and
\begin{equation*}
  S_2 =
  \sum_{n}
  \frac{\tilde{a}(n) a((n)_1)}{n (n)_1}
 \frac{(-1)^{\Omega(n)}}{\omega(n)}
 b(n) \log (n)_1.
\end{equation*}
Then we have
\begin{equation*}
\sum_{\substack{d \in \Omega_{-1}(a)}}
L'(1, E^{(d)})
E_{\ell}(-\mathcal{P}(d))
\Phi\left(\frac{-d}{X}\right)
\ll |S_1|X\log X+|S_2|X+X^{\frac{3}{4} + 2\varepsilon}.
\end{equation*}

For \( S_1 \), we use Rankin's trick to remove the restriction imposed by \( b(n) \).
When \( \Omega(n) > \ell \), we have \( 2^{\Omega(n) - \ell} > 1 \).
Thus, we obtain the bound
\begin{equation*}
  |S_1| \leq
  \left| \sum_{n}
  \frac{\tilde{a}(n) a((n)_1)}{n (n)_1}
 \frac{(-1)^{\Omega(n)}}{\omega(n)} \right|
  + \left| \sum_{n}
  \frac{|\tilde{a}(n) a((n)_1)| }{n (n)_1 }
  \frac{2^{\Omega(n) - \ell}}{\omega(n)}\right|.
\end{equation*}

For the first part, we have
\begin{multline*}
  \left| \sum_{n}
  \frac{\tilde{a}(n) a((n)_1)}{n (n)_1}
  \cdot \frac{(-1)^{\Omega(n)}}{\omega(n)} \right|
  =
  \left| \prod_{p \in \mathcal{P}}
  \left( \sum_{i=0}^{\infty}
  \frac{a(p)^{2i}}{p^{2i} (2i)!}
  - \sum_{i=0}^{\infty}
  \frac{a(p)^{2i+2}}{p^{2i+2} (2i+1)!} \right) \right|
  \\ \ll
   \prod_{p \in \mathcal{P}}
  \left( 1 - \frac{a(p)^2}{2p^2} + O\left( \frac{1}{p^2} \right) \right).
\end{multline*}

For the second part, using the inequality \( 1 + x \leq e^x \) for all \( x \in \mathbb{R} \), we obtain
\begin{multline*}
\left| \sum_{n}
  \frac{|\tilde{a}(n) a((n)_1)| }{n (n)_1 }
 \frac{2^{\Omega(n) - \ell}}{\omega(n)} \right|
  \leq
  2^{-\ell}
  \prod_{p \in \mathcal{P}}
  \left( \sum_{i=0}^{\infty}
  \frac{a(p)^{2i} 2^{2i}}{p^{2i} (2i)!}
  + \sum_{i=0}^{\infty}
  \frac{a(p)^{2i+2}  2^{2i+1}}{p^{2i+2} (2i+1)!} \right)
  \\ \leq
  2^{-\ell} \exp\left(
  2 \sum_{p \in \mathcal{P}}
  \left( \frac{a(p)^2}{p^2} + O\left( \frac{1}{p^2} \right) \right)
  \right).
\end{multline*}

Using Lemma~\ref{lemma:ortho}, we have
\begin{equation*}
 |S_1|\ll (\log X)^{-\frac{1}{2}}\log\log X.
\end{equation*}

By expanding the logarithmic term, we have
\begin{equation*}
S_2 =
 - \sum_{q \in \mathcal{P}}
 \sum_{l \geq 0}
 \frac{a(q)^{2l+2} \log q}{q^{2l+2} (2l+1)!}
 \sum_{\substack{n \\ (n,q) = 1}}
 \frac{\tilde{a}(n) a((n)_1)}{n (n)_1}
 \frac{(-1)^{\Omega(n)}}{\omega(n)}
b(n q^{2l+1}).
\end{equation*}
For the inner sum, when \( l \leq \frac{\ell}{2} \), we have
\begin{equation*}
   \left| \sum_{\substack{n \\ (n,q) = 1}}
   \frac{\tilde{a}(n) a((n)_1)}{n (n)_1}
\frac{(-1)^{\Omega(n)}}{\omega(n)}
   b(n q^{2l+1}) \right|
   \ll
   (\log X)^{-\frac{1}{2}} \log\log X
   \left( 1 + \frac{a(q)^2}{2q^2} + O\left( \frac{1}{q^2} \right) \right).
\end{equation*}

When \( l > \frac{\ell}{2} \), by applying Stirling’s formula to the factorial term, we obtain
\begin{multline*}
  \frac{1}{(2l+1)!} \left| \sum_{\substack{n \\ (n,q) = 1}}
  \frac{\tilde{a}(n) a((n)_1)}{n (n)_1}
  \frac{(-1)^{\Omega(n)}}{\omega(n)}
b(n q^{2l+1}) \right|
 \\ \ll
  \left( (\log X)^{-\frac{1}{2}} \log\log X +
  \frac{2^{2l+1}}{(2l+1)!} 2^{-\ell}
  (\log X)^2 (\log\log X)^{-4} \right)
  \\
  \times \left( 1 + \frac{a(q)^2}{2q^2} + O\left( \frac{1}{q^2} \right) \right)
  \\ \ll
  (\log X)^{-\frac{1}{2}} \log\log X
  \left( 1 + \frac{a(q)^2}{2q^2} + O\left( \frac{1}{q^2} \right) \right).
\end{multline*}

Then, applying Lemma~\ref{lemma:ortho}, we conclude that
\begin{equation*}
  |S_2| \ll
  \sum_{q \in \mathcal{P}}
  \left( \frac{a(q)^2}{q^2} + O\left( \frac{1}{q^2} \right) \right)
  (\log X)^{-\frac{1}{2}} \log\log X
  \ll
  (\log X)^{\frac{1}{2}} (\log\log X)^{-1}.
\end{equation*}

This completes the proof.
\end{proof}

\section*{Acknowledgements}

The author would like to thank Professor Shuai Zhai for pointing out the previous error.


\end{document}